\documentclass{elsart}

\usepackage[dvips]{epsfig}

\newcommand{\R}{I \! \! R}

\newcommand{\C}{I \! \! \! \! {C}}

\renewcommand{\o}{\omega}
\newcommand{\ba}{{\bf a}}
\newcommand{\by}{{\bf y}}
\newcommand{\bw}{{\bf w}}
\newcommand{\bvs}{{\bf {\underline a}}}
\newcommand{\bbf}{{\bf f}}

\newcommand{\bg}{{\mbox{\boldmath $\gamma$}}}
\newcommand{\bzet}{{\mbox{\boldmath $\zeta$}}}
\newcommand{\uu}{{\underline u}}
\newcommand{\ux}{{\underline x}}
\newcommand{\ut}{{\underline t}}

\newcommand{\uat}{\tilde{{\underline a}}}
\newcommand{\buat}{\tilde{\bf{\underline a}}}
\newcommand{\ue}{{\underline e}}
\newcommand{\vs}{{\underline a}}

\newcommand{\uzet}{{\underline{\zeta}}}
\newcommand{\ugam}{{\underline{\gamma}}}
\newcommand{\ugamt}{\tilde{{\underline{\gamma}}}}
\newcommand{\ugamh}{\hat{{\underline{\gamma}}}}

\newcommand{\us}{{\underline s}}
\newtheorem{teo}{Theorem}

\newcommand{\bb}{\begin{eqnarray*}}
\newcommand{\be}{\end{eqnarray*}}

\begin{document}

\markboth{P. Barone}
{Universality Pad\'{e} noise}

\title{On the universality of the distribution of the generalized eigenvalues of a pencil of Hankel random matrices}

\author{Piero Barone}

\address{ Istituto per le Applicazioni del Calcolo ''M. Picone'',
C.N.R.,\\
Via dei Taurini 19, 00185 Rome, Italy \\
barone@iac.rm.cnr.it, piero.barone@gmail.com}

\maketitle

\section*{Abstract}
Universality properties of the  distribution of the generalized eigenvalues of a pencil of random Hankel matrices, arising in the solution of the exponential interpolation problem of  a complex  discrete stationary process,  are proved under the  assumption that every finite set of random variables of the process have a multivariate spherical distribution. An integral representation of the condensed density of the generalized eigenvalues is also derived. The asymptotic behavior of this function turns out to depend only on stationarity and not on the specific distribution of the process.

{\it Key words}: complex moments; Pad\'{e}  approximants;  random polynomials.

{\it MSC 2000}: 15B52, 60B20, 62E15

\section*{Introduction}
Let us consider the following moment problem: to compute
 the complex measure defined on a compact set
$D\subset\C$ by
$$S(z)=\sum_{j=1}^p \gamma_j\delta(z-\zeta_j),\;\;\zeta_j\in \mbox{int}(
D), \;\;\zeta_j\ne\zeta_h \; \forall j\ne h,\;\;\gamma_j\in\C$$ from its complex moments
$$s_k=\int_Dz^kS(z)dz=\int\!\!\!\!\!\int_{D}(x+iy)^k S(x+iy)dx dy,
\;\;k=0,\dots,n-1,\,\;n=2p$$  It turns out that
\begin{eqnarray}s_k=\sum_{j=1}^p \gamma_j\zeta_j^k.\label{modale}\end{eqnarray}
The problem is thus equivalent to solve the complex exponential interpolation problem for the data $s_k,\;k=0,\dots,n-1$. It is well known that the $\zeta_j,\;j=1,\dots,p$ are the generalized eigenvalues of the Hankel pencil $P=[U_1(\us),U_0(\us)]$ where
 $$U_0(\us)=U(s_0,\dots,s_{n-2}),\;\;\;\;U_1(\us)=U(s_1,\dots,s_{n-1})$$
 and $$U(s_0,\dots,s_{n-2})=\left[\begin{array}{llll}
s_0 & s_{1} &\dots &s_{p-1} \\
s_{1} & s_{2} &\dots &s_{p} \\
. & . &\dots &. \\
s_{p-1} & s_{p} &\dots &s_{n-2}
  \end{array}\right].$$
Moreover $\gamma_j$ are related to the generalized
eigenvector $\uu_j$ of $P$ by $\gamma_j=\uu_j^T[s_0,\dots,s_{p-1}]^T$.

\noindent Equivalently, $\zeta_j,\;j=1,\dots,p$  are the roots of the polynomial  in
the variable $z$
$$det[U_1(\us)-zU_0(\us)]$$ which is the denominator of the Pad\'{e} approximant $[p-1,p](z)$ to the $Z$-transform of ${s_k}$ defined by $$f(z)=\sum_{k=0}^\infty s_kz^{-k}.$$

\noindent Denoting random variables by bold characters, let us assume now that we know an even number $n\ge 2p$ of noisy complex moments
$${\bf a}_k=s_k+\mbox{\boldmath $\nu$}_k,\quad k=0,1,2,\dots,n-1
$$ where $\mbox{\boldmath $\nu$}_k$  is a stationary discrete process, and we want to estimate $S(z)$ - i.e. $p,\;(\gamma_j,\zeta_j),j=1,\dots,p$ - from $\{{\bf
a}_k\}_{k=0,\dots,n-1}$. This is a well known difficult ill-posed problem which
is central in many disciplines and appears in the
literature in different forms and contexts (see e.g.
\cite{dsp,donoho,gmv,osb,scharf,vpb}). All the quantities defined above become random. It is therefore relevant to study the distribution of these quantities under suitable hypotheses on the noise affecting the data.

In this work the case $s_k=0,\;\forall k$ will be considered
assuming that the noise  is represented by a discrete stationary process, white or colored. For example the distribution of every finite set of r.v. of the process can be  multivariate  $\alpha$-stable which is a class of distributions closed with respect to addition (up to scale and location parameters), a property consistent with the naive concept of noise. Most of the properties proved in the following will require that the distribution of every finite set of r.v. of the process is spherical and the density function exists (\cite[Sec.1.5]{muir},\cite[Sec.2]{fang}). We remember that
under this assumption white noise is necessarily Gaussian. Therefore e.g. $\alpha$-stable, non-Gaussian $(\alpha\ne 2)$, spherical processes are colored.

Motivations to consider the pure noise case are twofold. The problem of identifying the presence of a signal in a large noise environment arises in many applied contexts. The  distribution of the generalized eigenvalues  $\mbox{\boldmath $\zeta$}_j$
of the random pencil
$${\bf
P}=[U({\bf a}_1,\dots,{\bf a}_{n-1}),U({\bf a}_0,\dots,{\bf
a}_{n-2})]$$ in the pure noise case is a reference to detect the presence of a signal.  More generally, when solving the noise filtering problem, information on the generalized spectral properties of the noise are required. In the specific class of complex exponential models a generalized spectrum is given by the condensed density (\cite{barjma}) of the generalized eigenvalues $\mbox{\boldmath $\zeta$}_j$
of the random pencil ${\bf P}$ or, equivalently, of the poles
 of the $[p-1,p](z)$ random Pad\'{e} approximant
to the $Z$-transform (formal random series)
$$\bbf(z)=\sum_{k=0}^\infty \ba_k z^{-k}.$$

In \cite{barja} it was proved that if $\{\ba_k\}$ is a complex Gaussian white noise, expressing the  condensed density in polar coordinates, the radial density weakly converges to a Dirac distribution centered in $1$ and the phase density is uniform in $(0,2\pi]$. Moreover it was proved that for $n=2$ the condensed  density is the uniform measure on the Riemann sphere
$$h_2(z)= \frac{1}{\pi(1+|z|^2)^2}.$$
In \cite{bess} it was conjectured, on the basis of numerical experiments, that the condensed density has an universal behavior i.e., in polar coordinates, the radial density is Lorentzian, centered in $1$ with width dependent on $n$ and the phase density is uniform in $(0,2\pi]$.

In the following we prove that this  conjecture is true under the hypothesis that the noise has a spherical distribution but the radial density is not Lorentzian. More specifically, assuming only the stationarity of the noise process,  it is proved that the condensed density is asymptotically concentrated on the unit circle independently of the noise distribution. Moreover it is proved that if the noise distribution is spherical
 the condensed density $\forall n$  is independent of the specific spherical distribution, it is invariant by scaling, and, when $n=2$, it is given by the uniform measure on the Riemann sphere. Furthermore in polar coordinates $(\rho,\theta)$ the marginal condensed density w.r. to $\rho$ (radial density) does not depend on $\theta$ and the marginal condensed density w.r. to $\theta$ is uniform
 .
 Finally an integral representation of the condensed density is provided in the spherical case which shows that  the  radial density is not Lorentzian.

The paper is organized in two sections. In the first one  the properties of the  condensed density are studied. In the second one three numerical experiments confirming the derived properties are illustrated.

\section{The condensed density of the Pad\'{e} poles}

Let us consider the transformation
$T=(T^{(1)},T^{(2)})$  that maps every
realization $\vs(\o)$ of $\bvs=\{\ba_k,k=0\dots,n-1\}$ to $(\uzet(\o), \ugam(\o))$
given by $
a_k(\o)=\sum_{j=1}^{n/2}\gamma_j(\o)\zeta_j(\o)^{k},\;\;k=0,\dots,n-1,$
where $\o\in\Omega$ and $\Omega$  is the space of events.
 It was proved in \cite[Lem.2]{barja2} that $T(\bvs)$ is defined and  one-to-one a.e.. By noticing that the $\mbox{\boldmath $\zeta$}_j$ are given by $T^{(1)}(\bvs)$,
 we define the  condensed density as
 $$h_n(z)=\frac{2}{n}E\left[\sum_{j=1}^{n/2}\delta(z-\bzet_j)\right]$$
 where the expectation is with respect to the density of $\bvs$.

\begin{teo}
If $\{\ba_k\}$ is stationary
$$\lim_{n \rightarrow\infty} h_n(z)=\delta(z-1)$$
\end{teo}

\noindent\underline{proof.}

Let us consider the periodic (circular) process $\{\tilde{\ba}_k\}$ obtained by repeating a finite segment of length $p$ of the process $\{\ba_k\}$. Then we have
$$\bbf(z)=\sum_{k=0}^\infty \tilde{\ba}_k z^{-k}= \sum_{h=0}^{p-1}\tilde{\ba}_h\left(\sum_{k=0}^\infty z^{-h+k p}\right)$$
if $|z|<1$ we get
$$\bbf(z)=\sum_{h=0}^{p-1}\tilde{\ba}_h z^{-h}\frac{z^p}{z^p-1}$$
hence $\bbf(z)$ is a random rational function for every period $p$ and its poles are the roots of unity. Therefore the Pade' denominator is a deterministic polynomial and  the poles' condensed density is the counting measure on its roots.
In the limit for $p\rightarrow \infty$ the periodic process $\{\tilde{\ba}_k\}$ converges in $L^2-$mean  to $\{\ba_k\}$ (\cite[Sec. 6]{kaw}) and  the counting measure on the roots of unity tends to the uniform measure on the unit circle.
$\Box$

\begin{teo}
If $n=2$ and $\buat=\{\Re[\ba_0],\Im[\ba_0],\Re[\ba_1],\Im[\ba_1]\}$ is spherically distributed with a density, then the poles'  condensed density is the uniform measure on the Riemann sphere $$h_2(z)=\frac{1}{\pi(1+|z|^2)^2}.$$

\end{teo}

\noindent\underline{proof}.

If a $n$-dimensional random vector $\buat$  is spherically distributed its  joint density is (\cite[th.2.9]{fang}) $$f(\|\uat\|^2)=\frac{\Gamma(n/2)}{2 \pi^{n/2}}\|\uat\|^{1-n}g(\|\uat\|)$$ where $g(\cdot)$ is the density of $\|\buat\|$ and $\|\cdot\|$ denotes the euclidean norm.  In the case considered $n=4$.
By making the change of variables $T$ we get $\bg=\ba_0,\;\;\bzet=\frac{\ba_1}{\ba_0}$ and
$$ \|\buat\|^2= |\bg|^2(1+|\bzet|^2).$$
Noticing that the complex Jacobian of $T$ is $\gamma$, the marginal on $\bzet$ is then
$$h(\zeta)=\int_{\C}\gamma f(T(\ux))d\gamma=\frac{1}{2 \pi^2}\int_{\C } \frac{\gamma}{|\gamma|^3(1+|\zeta|^2)^{3/2}}g\left(|\gamma|\sqrt{1+|\zeta|^2}\right)d\gamma.$$
By the change of variables $\tilde{\gamma}=\gamma\sqrt{1+|\zeta|^2)}$ and expressing the integral in real coordinates with real Jacobian $|\gamma|^2$ we get
$$h(\zeta)=\frac{1}{2 \pi^2}\int_{\R^2 }\frac{|\tilde{\gamma}|^2}{(1+|\zeta|^2)}\frac{1}{|\tilde{\gamma}|^3}g(|\tilde{\gamma}|)\frac{1}
{1+|\zeta|^2}d\Re{\tilde{\gamma}}d\Im{\tilde{\gamma}}=$$
$$\frac{1}{2 \pi^2(1+|\zeta|^2)^2}\int_{\R^2 }\frac{1}{|\tilde{\gamma}|}g(|\tilde{\gamma}|)d\Re{\tilde{\gamma}}d\Im{\tilde{\gamma}}.
$$
But $$\frac{1}{2\pi}\int_{\R^2 }\frac{1}{|\tilde{\gamma}|}g(|\tilde{\gamma}|)d\Re{\tilde{\gamma}}d\Im{\tilde{\gamma}}=1$$ because the marginals of a spherical density are spherical and we apply the formula above with $n=2\;\;\;\Box$

\begin{teo}
If $n> 2$ and $\buat=\{\Re[\bvs],\Im[\bvs]\}$ is $2n$-variate spherically distributed with a density, then the poles'  condensed density is the same independently of the specific distribution of $\buat$. Moreover in polar coordinates $(\rho,\theta)$ the marginal condensed density w.r. to $\rho$ does not depend on $\theta$ and the marginal condensed density w.r. to $\theta$ is uniform. Finally the condensed density is invariant by scaling $\buat$.

\end{teo}

\noindent\underline{proof}.

We first notice that the density of $\bvs$ is
\begin{eqnarray} f(\|\vs\|^2)=\frac{\Gamma(n)}{2 \pi^{n}}\|\vs\|^{1-2 n}g(\|\vs\|)\label{dens}\end{eqnarray} where $g(\cdot)$ is the density of $\|\bvs\|$ because
$\|\buat\|=\|\bvs\|$.
We have \bb
h_n(z)&=&\frac{2}{n}E\left[\sum_{j=1}^{n/2}\delta(z-\bzet_j)\right]\\&=&
\frac{\Gamma(n)}{2 \pi^{n}}\frac{2}{n}\sum_{j=1}^{n/2}\displaystyle
\int_{\C^n}\delta(z-\zeta_j)\|\vs\|^{1-2 n}g(\|\vs
\|)d\vs
\label{dc1}\be
and, by making the change of variables $T$, whose complex Jacobian is (\cite[Th.2]{barja2})
$$ J_C(\uzet,\ugam)=
 (-1)^{n/2}\prod_{j=1}^{n/2}\gamma_j
\prod_{j<h}(\zeta_j-\zeta_h)^4, $$
 we get
\bb
h_n(z)&=&
\frac{2}{n}\sum_{j=1}^{n/2}\frac{\Gamma(n)}{2\pi^{n}}\displaystyle
\int_{\C^{n/2}}\int_{\C^{n/2}}\delta(z-\zeta_j)J_C(\uzet,\ugam)\left(\sum_{k=0}^{n-1}\left|
\sum_{h=1}^{n/2}\gamma_h\zeta_h^{k}\right|^2\right)^{1/2- n}\cdot \\ &&g\left(\left(\sum_{k=0}^{n-1}\left|
\sum_{h=1}^{n/2}\gamma_h\zeta_h^{k}\right|^2\right)^{1/2}\right)d\uzet
d\ugam\\
&=&
\frac{2}{n}\sum_{j=1}^{n/2}\frac{\Gamma(n)}{2\pi^{n}}
\int_{\C^{n/2-1}}\int_{\C^{n/2}}
J_C^*(\uzet^{(j)},z,\ugam)
\left(\sum_{k=0}^{n-1}\left|
\sum_{h\ne j}^{1,n/2} \gamma_h\zeta_h^{k}+\gamma_j
z^{k}\right|^2\right)^{1/2- n}\cdot \\ &&g\left(\left(\sum_{k=0}^{n-1}\left|
\sum_{h\ne j}^{1,n/2} \gamma_h\zeta_h^{k}+\gamma_j
z^{k}\right|^2\right)^{1/2}\right)
d\uzet^{(j)}d\ugam
\label{dc2}\be
where $\uzet^{(j)}=\{\zeta_h,h\ne j\}$
and
$$J_C^*(\uzet^{(j)},z,\ugam)=
(-1)^{n/2}\prod_{ h=1}^{1,n/2}\gamma_h\prod_{r<h,r\ne
j}(\zeta_r-\zeta_h)^4\prod_{(r,h)\ne j}(\zeta_r-z)^4 .$$
Let us define
$$Q_j=Q_j(\uzet^{(j)},z)=X_jX^H_j\in\C^{n/2\times n/2},\;\;X_j\in\C^{n/2\times n},\;\;X_j(h,k)=\overline{x}_{hk}^{(j)}$$ where
\begin{eqnarray}x_{hk}^{(j)}=\left\{\begin{array}{ll}
\zeta_h^{k-1}, & h\ne j\\
z^{k-1}, & h=j
\end{array}\right. .\label{matx}\end{eqnarray}
Then
$$\sum_{k=0}^{n-1}\left|
\sum_{h\ne j}^{1,n/2} \gamma_h\zeta_h^{k}+\gamma_j
z^{k}\right|^2=\ugam^HQ_j\ugam$$
and
\begin{equation}
h_n(z)=\frac{2}{n}\sum_{j=1}^{n/2}h^{(j)}_n(z)\label{dc}\end{equation}
where
\begin{eqnarray} h^{(j)}_n(z)=
\int_{\C^{n/2-1}}\prod_{r<h,(r,h)\ne
j}(\zeta_r-\zeta_h)^4\prod_{r\ne j}(\zeta_r-z)^4G_j(\uzet^{(j)},z) d\uzet^{(j)}
\label{condh}\end{eqnarray} and
$$G_j(\uzet^{(j)},z)=\frac{\Gamma(n)}{2 \pi^{n}}\int_{\C^{n/2}}
(-1)^{n/2}\left(\prod_{ h=1}^{1,n/2}\gamma_h\right)
\left(\ugam^HQ_j\ugam\right)^{1/2- n}g\left(\left(\ugam^HQ_j\ugam\right)^{1/2}\right)
d\ugam.$$
We show now that $G_j(\uzet^{(j)},z)$ is independent of the specific form of the spherical density generator $f(\cdot)$.
If $\ugamt=\{\Re[\ugam],\Im[\ugam]\}\in \R^{n}$ and if $\tilde{Q}_j$ is the real isomorph of $Q_j$, i.e.
$$\tilde{Q}_j=\left[\begin{array}{ll}
\Re[Q_j], & -\Im[Q_j]\\
\Im[Q_j], & \;\;\;\Re[Q_j]
\end{array}\right],$$
 we have
$\ugam^H Q_j\ugam=\ugamt^T\tilde{Q}_j\ugamt$. But then, by noticing that  the real Jacobian of the transformation $T$  is $J_R=|J_C|^2$,
we get
$$G_j(\uzet^{(j)},z)=\frac{\Gamma(n)}{2 \pi^{n}}\int_{\R^n}
\left(\prod_{ h=1}^{1,n/2}\ugamt^T A_h\ugamt\right)
\left(\ugamt^T\tilde{Q}_j\ugamt\right)^{1/2- n}g\left(\left(\ugamt^T\tilde{Q}_j\ugamt\right)^{1/2}\right)
d\ugamt$$
where $A_h=I_2\otimes\ue_h\ue_h^T$ and $\ue_h$ is the $h-$ column of the identity matrix of order $n/2$.
Let us consider the transformation
$$\ugamh=\tilde{Q}_j^{1/2}\ugamt$$
which is well defined because $\tilde{Q}_j>0$ and whose Jacobian is $|\tilde{Q}_j|^{1/2}$. We then have
$$G_j(\uzet^{(j)},z)=\frac{\Gamma(n)}{2 \pi^{n}|\tilde{Q}_j|^{1/2}}\int_{\R^n}
\left(\prod_{ h=1}^{1,n/2}\ugamh^T\tilde{Q}_j^{-1/2} A_h\tilde{Q}_j^{-1/2}\ugamh\right)
(\ugamh^T\ugamh)^{1/2- n}g\left(\left(\ugamh^T\ugamh\right)^{1/2}\right)
d\ugamh=$$
$$\frac{\Gamma(n)}{\Gamma(n/2) \pi^{n/2}|\tilde{Q}_j|^{1/2}}\frac{\Gamma(n/2)}{2 \pi^{n/2}}\int_{\R^n}
\left(\prod_{ h=1}^{1,n/2}\frac{\ugamh^T\tilde{Q}_j^{-1/2} A_h\tilde{Q}_j^{-1/2}\ugamh}{\ugamh^T\ugamh}\right)
\left(\sqrt{\ugamh^T\ugamh}\right)^{1- n}g\left(\left(\ugamh^T\ugamh\right)^{1/2}\right)
d\ugamh.$$
But then
$$\frac{\Gamma(n/2)}{2 \pi^{n/2}}\int_{\R^n}
\left(\prod_{ h=1}^{1,n/2}\frac{\ugamh^T\tilde{Q}_j^{-1/2} A_h\tilde{Q}_j^{-1/2}\ugamh}{\ugamh^T\ugamh}\right)
\left(\sqrt{\ugamh^T\ugamh}\right)^{1- n}g\left(\left(\ugamh^T\ugamh\right)^{1/2}\right)
d\ugamh$$ can be seen as the expectation of the function $\left(\prod_{ h=1}^{1,n/2}\frac{\ugamh^T\tilde{Q}_j^{-1/2} A_h\tilde{Q}_j^{-1/2}\ugamh}{\ugamh^T\ugamh}\right)$ of a $n-$dimensional  spherically distributed r.v.
$\bg$. Therefore
$$G_j(\uzet^{(j)},z)=\frac{\Gamma(n)}{\Gamma(n/2) \pi^{n/2}|\tilde{Q}_j|^{1/2}}E\left[\prod_{ h=1}^{1,n/2}\left(\frac{\bg^T\tilde{Q}_j^{-1/2} A_h\tilde{Q}_j^{-1/2}\bg}{\bg^T\bg}\right)\right].$$
Let us define the  matrix
$$B_{jh}=\frac{\tilde{Q}_j^{-1/2} A_h\tilde{Q}_j^{-1/2}}{\Re[Q_j^{-1}]_{hh}}$$
and prove that $B_{jh}$ is idempotent. In fact
$$B_{jh}^2=\frac{\tilde{Q}_j^{-1/2} A_h\tilde{Q}_j^{-1/2}}{\Re[Q_j^{-1}]_{hh}}\cdot \frac{\tilde{Q}_j^{-1/2} A_h\tilde{Q}_j^{-1/2}}{\Re[Q_j^{-1}]_{hh}}=\frac{\tilde{Q}_j^{-1/2} A_h\tilde{Q}_j^{-1}A_h\tilde{Q}_j^{-1/2}}{\Re[Q_j^{-1}]_{hh}^2}=$$ $$
\frac{\tilde{Q}_j^{-1/2}(I_2\otimes\ue_h\ue_h^T)\tilde{Q}_j^{-1}(I_2\otimes\ue_h\ue_h^T)\tilde{Q}_j^{-1/2}}{\Re[Q_j^{-1}]_{hh}^{2}}.$$
But
\bb(I_2\otimes\ue_h\ue_h^T)\tilde{Q}_j^{-1}(I_2\otimes\ue_h\ue_h^T)&=&\left[\begin{array}{ll}
\ue_h\ue_h^T, & 0\\
0, & \ue_h\ue_h^T
\end{array}\right]
\left[\begin{array}{ll}
\Re[Q_j^{-1}], & -\Im[Q_j^{-1}]\\
\Im[Q_j^{-1}], & \;\;\;\Re[Q_j^{-1}]
\end{array}\right]\left[\begin{array}{ll}
\ue_h\ue_h^T, & 0\\
0, & \ue_h\ue_h^T
\end{array}\right]\\&=& \Re[Q_j^{-1}]_{hh}A_h\be because the diagonal elements of $\Im[Q_j^{-1}]$ are zero.
Hence
$$B_{jh}^2=\tilde{Q}_j^{-1/2}\frac{\Re[Q_j^{-1}]_{hh}}{(\Re[Q_j^{-1}]_{hh})^{2}}A_h\tilde{Q}_j^{-1/2}=B_{jh}.$$
We then have
\begin{eqnarray}G_j(\uzet^{(j)},z)=\frac{\Gamma(n)}{\Gamma(n/2) \pi^{n/2}|\tilde{Q}_j|^{1/2}}\prod_{ h=1}^{1,n/2}[Q_j^{-1}]_{hh}E\left[\prod_{ h=1}^{1,n/2}\left(\frac{\bg^TB_{jh}\bg}{\bg^T\bg}\right)\right].\label{condG}\end{eqnarray}
From \cite[Th.1.5.7,ii]{muir} it follows that $\bw_h=\left(\frac{\bg^T B_{jh}\bg}{\bg^T\bg}\right)$ has the beta distribution with parameters $1$ and $n/2-1$ independently of the distribution of $\bg$. As a consequence the distribution of  $h_n(z)$ too does not depend on the distribution of $\bg$.

\noindent To prove the last part of the theorem, we notice that the spherical density (\ref{dens}) is invariant under the transformation
$$\vs\rightarrow e^{\pm i \frac{\beta}{2}}\vs,\;\;\forall \beta.$$
The proof of Theorem 2 in \cite{barja} then holds and provides the requested results. $\;\;\;\Box$

\begin{teo}
If $n> 2$ and $\buat=\{\Re[\bvs],\Im[\bvs]\}$  is $2n$-variate spherically distributed with a density, the poles condensed density is given by
\bb h_n(z)=\frac{1}{(2\pi)^{n/2}}
\int_{\R^{n-2}}K_{G_n}(z,\uzet^{(1)})\prod_{r<h,(r,h)\ne
1}|\zeta_r-\zeta_h|^{2}\prod_{r\ne 1}|\zeta_r-z|^{2}\cdot\\
\frac{\prod_{ h=1}^{1,n/2}\left(\sum_{P_2}\left|s_{P_2}(z,\uzet^{(1,h)})\right|^2\right)}{ \left(\sum_{P_1}\left|s_{P_1}(z,\uzet^{(1)})\right|^2\right)^{n/2+1}}d\Re(\uzet^{(1)})d\Im(\uzet^{(1)})\be
where
$$K_{G_n}(z,\uzet^{(1)})=E\left[\prod_{ h=1}^{1,n/2}\left(\by^TB_{1h}\by\right) \right]$$
and $\by$ is a $n$-variate standard Gaussian random vector;

\noindent $s_{P_1}(z,\uzet^{(1)})$ are the Schur functions associated to the partition $$P_1=\{j_1 , j_2 , \dots , j_{n/2}\}$$ which spans the minors of maximal order  of $X_1$ (eq. \ref{matx});

\noindent
$s_{P_2}(z,\uzet^{(1,h)})$ are the Schur functions associated to the partition $$P_2=\{j_1 , j_2 , \dots , j_{n/2-1}\}$$ which spans the minors of maximal order  of the matrix obtained by $X_1$ by canceling the $h-$th row.
\end{teo}

\noindent\underline{proof}.

From   Theorem 3, without loss of generality, we can assume that $\bg$ has a zero-mean multivariate Gaussian density with identical covariance matrix. The ratio of quadratic forms $\bw_h=\left(\frac{\bg^TB_{jh}\bg}{\bg^T\bg}\right)$ can be rewritten as  $\bw_h=\by^TB_{jh}\by$ where $\by=\frac{\bg}{\|\bg\|}$. Therefore $\bw_h$ is a quadratic form in the variables $\by$ which are uniformly distributed on the $n-$dimensional sphere. But this distribution is spherical with characteristic function
given by (\cite[3.1.1]{fang}
$$\Psi_n(\ut)=\, _0F_1\left(\frac{n}{2};\frac{1}{4}\|\ut\|^2\right).$$
In \cite[Prop.4]{kan} an explicit expression for
$K_{G_n}(z,\uzet^{(j)})=E\left[\prod_{ h=1}^{1,n/2}\bw_h\right]$ when $\by$ is multivariate Gaussian is given, and it is claimed that this is also the result in the spherical case up to a constant which is a function of the characteristic function of the spherical distribution. It turns out that in the present case the constant is given by
$$\beta_n=\frac{n^{n/2-1}}{4}\frac{\Psi_n^{(n/2+1)}(0)}{\left(\Psi_n^{(1)}(0)\right)^{n/2}}=\frac{\Gamma[n/2]}{2^{n/2}\Gamma[n]}.$$
 We then have
$$E\left[\prod_{ h=1}^{1,n/2}\bw_h\right]= \beta_n\cdot K_{G_n}(z,\uzet^{(j)})$$
 Moreover this expression is a symmetric function of $z,\uzet^{(j)}$.
We then get by  equations (\ref{condh}) and (\ref{condG})
\bb h^{(j)}_n(z)&=&\frac{1}{(2\pi)^{n/2}}
\int_{\C^{n/2-1}}K_{G_n}(z,\uzet^{(j)})\prod_{r<h,(r,h)\ne
j}(\zeta_r-\zeta_h)^4\prod_{r\ne j}(\zeta_r-z)^4\cdot \\ &&\frac{1}{|\tilde{Q}_j|^{1/2}}
\prod_{ h=1}^{1,n/2}[Q_j^{-1}]_{hh} d\uzet^{(j)}.\be
Remembering  that the real Jacobian of the transformation $T$ is $J_R=|J_C|^2$, we have
\bb h^{(j)}_n(z)=\frac{1}{(2\pi)^{n/2}}
\int_{\R^{n-2}}K_{G_n}(z,\uzet^{(j)})\prod_{r<h,(r,h)\ne
j}|\zeta_r-\zeta_h|^8\prod_{r\ne j}|\zeta_r-z|^8\frac{1}{|Q_j|}\cdot \\
\prod_{ h=1}^{1,n/2}[Q_j^{-1}]_{hh} d\Re(\uzet^{(j)})d\Im(\uzet^{(j)})\be
because  $|\tilde{Q_j}|=|Q_j|^2$.
 But, by Binet-Cauchy formula, we have
$$|Q_j|=\sum_{j_1\le j_2<\dots<j_{n/2}}\left|X_j\left(\begin{array}{llll}1 & 2 & \dots & n/2 \\ j_1 & j_2 & \dots & j_{n/2}\end{array}\right)\right|^2$$
where $X_j\left(\begin{array}{llll}1 & 2 & \dots & n/2 \\ j_1 & j_2 & \dots & j_{n/2}\end{array}\right)$ is a minor of maximal order $n/2$ of $X_j$. From \cite{mcdon} we have
$$X_j\left(\begin{array}{llll}1 & 2 & \dots & n/2 \\ j_1 & j_2 & \dots & j_{n/2}\end{array}\right)=s_{P_1}(z,\uzet^{(j)})\prod_{r<h,(r,h)\ne j}(\zeta_r-\zeta_h)\prod_{r\ne j}(\zeta_r-z)$$
where $s_{P_1}(z,\uzet^{(j)})$ is the Schur function associated to the partition $P_1=\{j_1 , j_2 , \dots , j_{n/2}\}$, which is a symmetric polynomial with positive integer coefficients (Jack function with $\alpha=1$).
Hence
$$|Q_j|=\prod_{r<h,(r,h)\ne j}|\zeta_r-\zeta_h|^2\prod_{r\ne j}|\zeta_r-z|^2 \sum_{j_1\le j_2<\dots<j_{n/2}}\left|s_{P_1}(z,\uzet^{(j)})\right|^2.$$
 But
 $[Q^{-1}_j]_{hh}=\frac{[\mbox{adj}( Q_j)]_{hh}}{|Q_j|}$ and
 $$[\mbox{adj}( Q_j)]_{hh}=\sum_{j_1\le j_2<\dots<j_{n/2-1}}\left|X_j\left(\begin{array}{lllllll}1, & 2, &\dots & h-1,&h+1\dots & n/2 \\ j_1, & j_2, & \dots& \dots& \dots & j_{n/2-1}\end{array}\right)\right|^2=$$
 $$\sum_{P_2}\left|s_{P_2}(z,\uzet^{(j,h)})\right|^2\prod_{r<k,(r,k)\ne j,h}|\zeta_r-\zeta_k|^2\prod_{r\ne j\ne h}|\zeta_r-z|^2$$
 where $P_2=\{j_1 , j_2 , \dots , j_{n/2-1}\}$
 and
\bb h^{(j)}_n(z)= \frac{1}{(2\pi)^{n/2}}
\int_{\R^{n-2}}K_{G_n}(z,\uzet^{(j)})\prod_{r<h,(r,h)\ne
j}|\zeta_r-\zeta_h|^{6-n}\prod_{r\ne j}|\zeta_r-z|^{6-n}\cdot\\
\prod_{ h=1}^{1,n/2}\left(\prod_{r<k,(r,k)\ne j,h}|\zeta_r-\zeta_k|^2\prod_{r\ne j\ne h}|\zeta_r-z|^2 \right)\cdot\\
\frac{\prod_{ h=1}^{1,n/2}\left(\sum_{P_2}\left|s_{P_2}(z,\uzet^{(j,h)})\right|^2\right)}{ \left(\sum_{P_1}\left|s_{P_1}(z,\uzet^{(1)})\right|^2\right)^{n/2+1}}d\Re(\uzet^{(j)})d\Im(\uzet^{(j)}).\be
But it turns out that
$$\prod_{ h=1}^{1,n/2}\left(\prod_{r<k,(r,k)\ne j, h}|\zeta_r-\zeta_k|^2\prod_{r\ne j\ne h}|\zeta_r-z|^2 \right)=\prod_{r<h,(r,h)\ne
j}|\zeta_r-\zeta_h|^{n-4}\prod_{r\ne j}|\zeta_r-z|^{n-4}$$ therefore
\bb h^{(j)}_n(z)= \frac{1}{(2\pi)^{n/2}}
\int_{\R^{n-2}}K_{G_n}(z,\uzet^{(j)})\prod_{r<h,(r,h)\ne
j}|\zeta_r-\zeta_h|^{2}\prod_{r\ne j}|\zeta_r-z|^{2}\cdot\\
\frac{\prod_{ h=1}^{1,n/2}\left(\sum_{P_2}\left|s_{P_2}(z,\uzet^{(j,h)})\right|^2\right)}{ \left(\sum_{P_1}\left|s_{P_1}(z,\uzet^{(1)})\right|^2\right)^{n/2+1}}d\Re(\uzet^{(j)})d\Im(\uzet^{(j)}).\be
Because of the symmetry of the Schur polynomials and of $K_n(z,\uzet^{(j)})$ all $h^{(j)}_n(z),\;j=1,\dots,n/2$ are equal. Therefore by equation (\ref{dc}) we get the thesis.
$\;\;\;\Box$

From the  result of the theorem it follows that, expressing $h_n(z)$ in polar coordinates and taking the marginal with respect to the modulus (radial density), we get a much more complicated expression than a Lorentzian function, thus disproving a part of the conjecture mentioned in the Introduction. We also notice that evaluating $h_n(z)$ is not an easy task because computing Schur functions is far from trivial (\cite{sch}).

\section{Numerical examples}
 In order to illustrate the results obtained in Theorem 3 and 4 three  numerical experiments were performed. In the first one a sample of cardinality $2\cdot 10^6$ of a multivariate complex $\alpha-$stable, centered, symmetric colored noise with $\alpha=0.5$ and $n=4$ was generated. The Pade' poles were computed as well as the empirical density of their modulus. In fig.1 this estimate is represented by dots. Then the integral representation of the poles condensed density derived in Theorem 4 was transformed in polar coordinates $z=\rho\cos(\theta)$ with Jacobian $\rho$. After Theorem 3 we know
  that the phase distribution is uniform. Therefore the radial distribution is obtained by multiplying by $2\pi\rho$ the integral representation given in Theorem 4. The integral was approximated by numerical quadrature where, in this case $$K_{G_4}(z,\uzet^{(j)})= \mbox{tr}(B_1)\mbox{tr}(B_2)+2\mbox{tr}(B_1B_2),$$ and is plotted in fig.1 as a solid line. We notice that the fit is quite accurate.

  In the second experiment a $4$-variate complex Gaussian white noise was generated. The same computations as above were performed and the results are plotted in fig.2 confirming that the condensed density in the Gaussian case is the same than in the $\alpha-$stable one as claimed in Theorem 3.

  In the third experiment a sample of cardinality $2\cdot 10^6$ of a $6$-variate complex Gaussian white noise was generated.
  The same computations as above were performed. In this case
  \bb K_{G_6}(z,\uzet^{(j)})=\mbox{tr}(B_1)\mbox{tr}(B_2)\mbox{tr}(B_3)+2\left[\mbox{tr}(B_1B_2)\mbox{tr}(B_3)\right.+\\\left.
  \mbox{tr}(B_1B_3)\mbox{tr}(B_2)+\mbox{tr}(B_2B_3)\mbox{tr}(B_1)\right]+8\mbox{tr}(B_1B_2B_3).\be The results are plotted in fig.3. Also in this case the fit is quite accurate.


\begin{figure}[H]
\hspace*{-0.5in}
\includegraphics[totalheight=5.in]{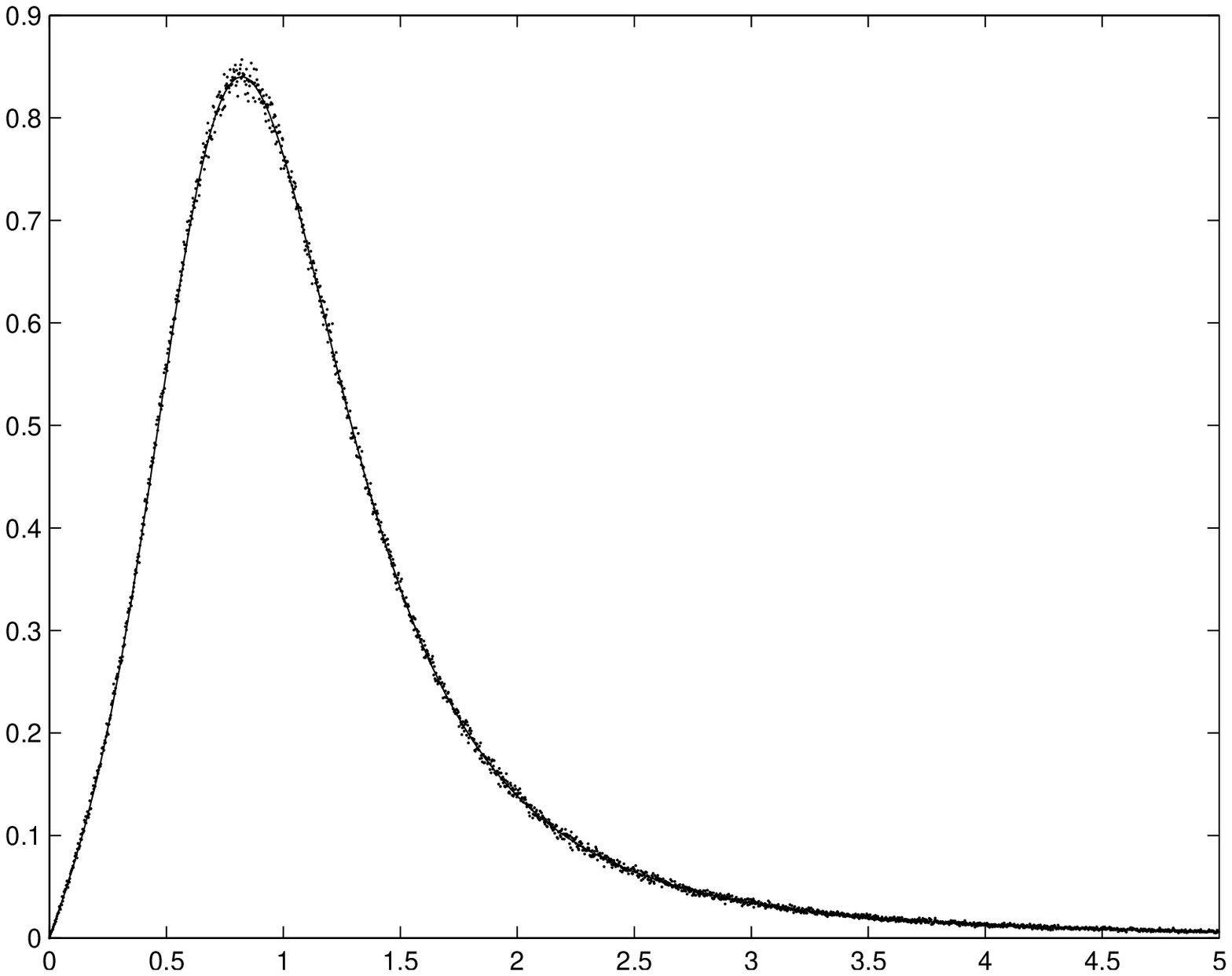}
\caption{The poles condensed density for $n=4$, computed by  numerical integration (Th.4) (solid) and the empirical density based on $2\cdot 10^6$ replications of a colored noise with a symmetric, centered $\alpha-$stable density with $\alpha=0.5$ (dotted).}
\label{fig1}
\end{figure}


\begin{figure}[H]
\hspace*{-0.5in}
\includegraphics[totalheight=5.in]{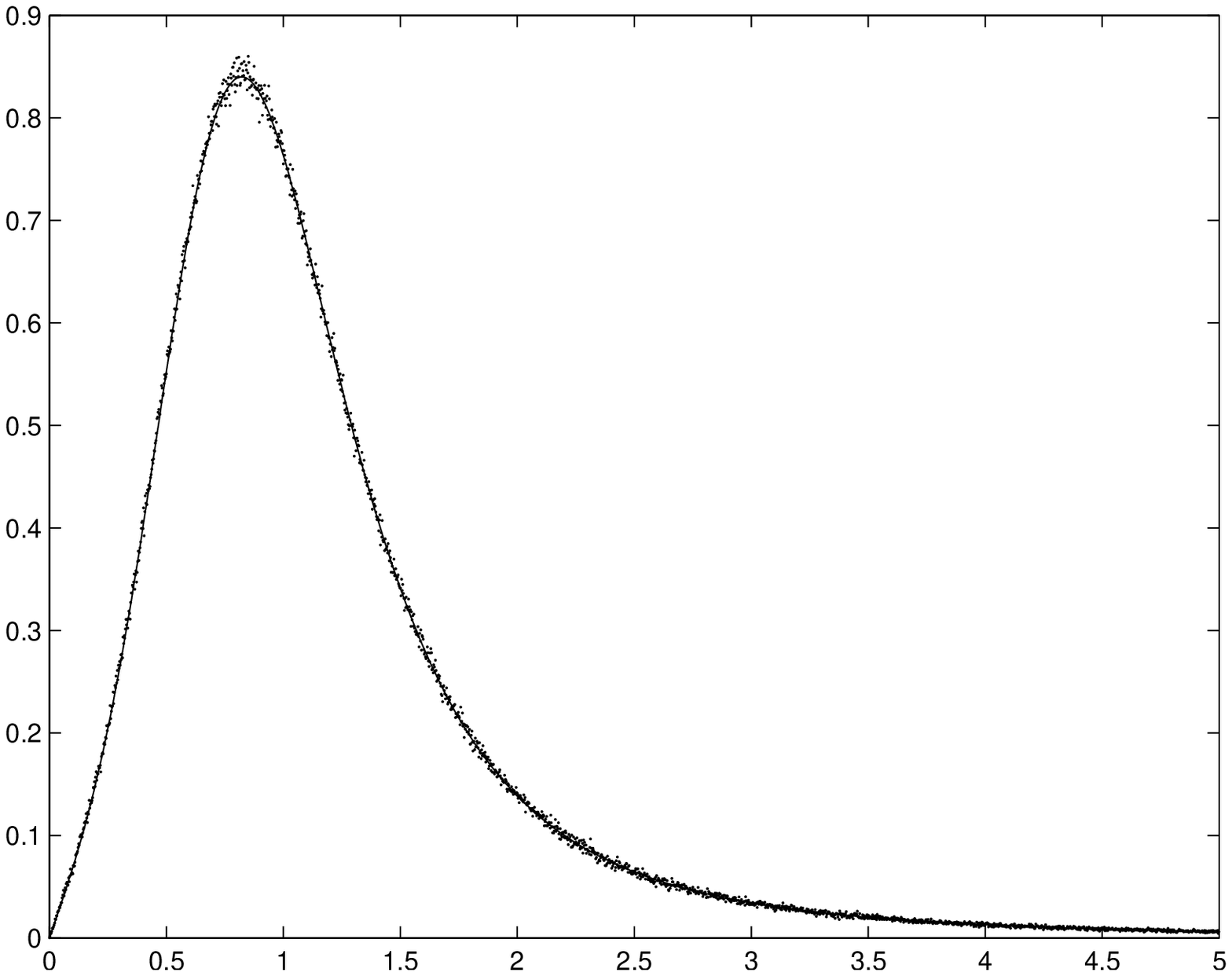}
\caption{The poles condensed density for $n=4$, computed by  numerical integration (Th.4) (solid) and the empirical density based on $2\cdot 10^6$ replications of a complex Gaussian white noise (dotted).}
\label{fig2}
\end{figure}


\begin{figure}[H]
\hspace*{-0.5in}
\includegraphics[totalheight=5.in]{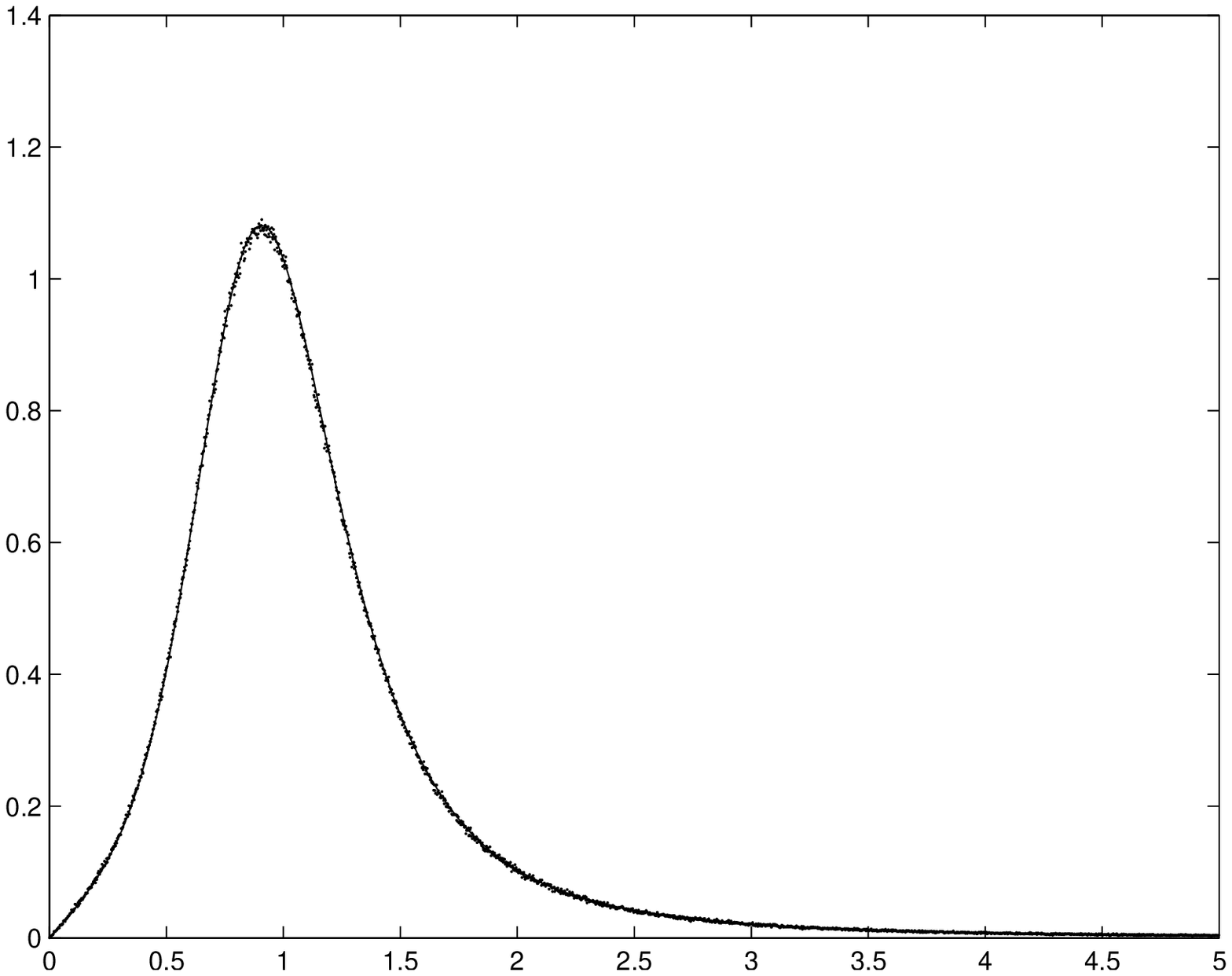}
\caption{The poles condensed density for $n=6$, computed by  numerical integration (Th.4) (solid) and the empirical density based on $2\cdot 10^6$ replications of a complex Gaussian white noise (dotted).}
\label{fig3}
\end{figure}

\end{document}